\documentclass{article}

\usepackage{amsmath}
\usepackage{graphicx}
\usepackage{subcaption}
\usepackage[a4paper, total={6in, 8in}]{geometry}
\usepackage{tabularx}
\usepackage{upgreek}
\usepackage{amssymb}
\usepackage{authblk}

\begin{document} 

\title{Stabilizing effect of delay in higher dimensions} 
\author{Alena Chan}
\affil{alenachan121@gmail.com}
\date{September 1, 2022} 
\maketitle 
 
    \section{Introduction}
    
    \text In this paper we study the stabilization effect of delay in dynamical system with three agents. The idea of expanding the system in [1] is motivated by desire to model more complex interactions among many characters. Firstly, we can model relationships between three agents, following Strogatz's [2] idea of 'many-body problem'. We consider a situation where changes of love/hate in time are described as a system of three linear differential equations, and assume that the system 
    has an unstable equilibrium. Our goal is to make the system stable by 
    introducing an adjustment: at least one of the agents responds not instantly 
    but with a certain delay. Thus, we replace the original system of equations
    by a system of \textit{delayed} differential equations. 
    In addition to the situation just described, delayed models can represent population dynamics, traffic control, red blood cell circulation, payload oscillations, stock markets, thermochemical reactions, etc. - see more detailed examples in [5]. \\
    
    \noindent
    The aim of this paper is to find conditions required for three-variable systems to gain stability. \\

    \noindent
    We now briefly outline the main ideas in our approach. 
    As can be seen from [4], Theorem 4.5, the system will be stable if its characteristic function $W(\lambda)$ has no zeros in the right half plane. 
    In general, the number $N$ of roots in the right half of complex plane is bounded, and we want $N$ to be $0$; as for the left half, there may be an infinite number of roots, but they do not affect stability. Hence, our aim is to estimate the number of roots with positive real part. We start with $\tau=0$, assuming that the steady state is unstable with $N=2$ roots in right half plane. \\
    
    \noindent 
    From the Argument Principle, $N$ can be evaluated by computing $\Updelta_{C} \arg{W(i \omega)}$ on the boundary $C$ of the right half circle with center at origin and radius $R >> 0$ and multiplying the value by $-\frac{1}{2\pi}$. Note that as the contour is directed clockwise, the factor is negative. 
    By arguments similar to those used in the proof of Rouche's Theorem, the change of argument on the curvilinear part is $-3\pi$. Moreover, the changes of two vertical half lines are equal since the characteristic function $W(\lambda)$
    is assumed to have real coefficients. \\
    
    \noindent 
    Hence we want the change of argument for the curve $W(i \omega)$ for $\omega \in [0; + \infty)$ to be $\frac{3 \pi}{2}$. As $\tau$ changes from $0$ to $+\infty$ the curve $W(i \omega)$ changes smoothly and $\Updelta \arg{W(i \omega)}$ can only change when $\tau$ goes through a value admitting 
    a solution to $W(i \omega)=0$ for some $\omega$. Then the argument changes by a multiple of $2 \pi$. When $\tau=0$, the change of argument is $-\frac{\pi}{2}$,
    and we want to know whether for a positive value of $\tau$ it can become $ \frac{3\pi}{2}$.
    Generalizing the approach of [1], p.3930 we will show that the values of $\tau$
    for which $W(i\omega) = 0$ admits a solution, come from two arithmetic series
    $$
    \tau_{1n} = \tau_{10} + \frac{2\pi n}{\omega_1}, \qquad
    \tau_{2n} = \tau_{20} + \frac{2\pi n}{\omega_2}, \qquad n \geq 0, 
    \qquad 0 < w_1 < w_2
    $$
    As the value of $\tau$ goes through $\tau_{1n}$ the number of roots on the right goes from $N$ to $N-2$ (situation gets 'better') and when $\tau$ goes through $\tau_{2n}$ the number of roots in the the right half plane goes from $N$ to $N+2$ (situation gets 'worse'). Note that we choose coefficients so that for $\omega=0$, $W(0) \neq 0$. Also, as $\omega_2 > \omega_1$, the situation gets 'worse' more often than it
    gets 'better' and eventually the system becomes unstable for large values of the delay $\tau$.\\
    
    \noindent
    In the case when $\tau_{10} < \tau_{20}$ the system goes from the unstable
    equilibrium ($N=2$) to stable ($N=0$), then again to unstable and perhaps
    after a few further alterations gets stuck in the unstable mode. When 
    $\tau_{10} > \tau_{20}$, the system goes from $N=2$ to $N=4$ and getting
    'better' always occurs after getting 'worse' hence the stable equilibrium never
    happens at all. \\

    \noindent 
    Our main result, Theorem, presents a series of conditions on the coefficients
    that guarantee the first of two scenarios. From the above, it suffices 
    to find any value $\overline{\tau}$ for which we have stability (not necessarily
    within the initial stability interval $[\tau_{10}, \tau_{20}]$). For that 
    we will formulate some conditions that ensure that the curve $W(i\omega)$
    starts at a point on the positive half of the horizontal axis, then moves through Quadrants I, then II, and continues to
    Quadrant III - receiving a contribution of $\frac{\pi}{2}$ from each of the quadrants. \\
    
     \noindent
    This paper is organized as follows. 
    In Section 2 we introduce a delayed model and describe the meaning of variables. Then in Section 3 we prove Theorem that guarantees occurrence of stability switches. We conclude that values of delay within a particular range result in stabilization. 
    Next, we summarize all conditions required for stabilization. In Section 4 we provide an example satisfying conditions of Theorem and show visualizations of stability switches. We conclude the paper by suggesting further open questions in Section 5.\\
    
    \noindent 
    \textbf{Acknowledgements.} The author thanks Prof. A. Gorodetski for formulating
    the problem and pointing to a few valuable references and Prof. V. Baranovsky 
    for helping to master the required background and for his continued support throughout the different stages of the project.

    \section{Romeo, Paris and Juliet model}
    
    \text We follow the ideas of [1] by expanding the model from two to three characters.
    We assume that for $\tau=0$ the system is unstable, which we wish to fix by introducing a positive delay. \\
    
    \noindent
    \textbf{Delayed model.} Our three-agent model reads
    \begin{equation}
        \begin{cases}
          x'(t) = -a_{0}z(t)-b_{0}z(t-\tau) + A_{1} \\
          y'(t) = -x(t) + a_{1}z(t) + b_{1}z(t-\tau)+ A_{2} \\
          z'(t) = -y(t) - a_{2}z(t) - b_{2}z(t-\tau) + A_{3}
        \end{cases}\,
    \end{equation}
    \text where $x(t)$ denotes Romeo’s emotions (love if $x(t)>0$, hate if $x(t)<0$) for Juliet at time $t$, $y(t)$ denotes Paris’ love/hate
    for Juliet at time $t$ and $z(t)$ denotes Juliet’s feelings for Romeo at time $t$. $A_{i}$ are appeal terms which are constant. Note that initial conditions for delayed differential equations should be defined in the whole interval $t \in [0; \tau]$; however, as the proof of Theorem 4.5 in [4] does not use initial conditions, they won't affect the stability.
    
    \section{General theorem}
    \text We prove that in the case described in Eqs. (1), under certain conditions on the coefficients there are two critical values of the delay $\tau$: the first stabilizes and the second destabilizes the system. The derivation of  quasi-polynomial characteristic equation can be found in Appendix, where we write solutions in the form $e^{\lambda t} * $ (constant vector). We assume that
    for our system of three delayed differential equations with one delay its characteristic function has a  form
    \begin{equation}
        W(\lambda) = \lambda^3 + a_{2}\lambda^2 + a_{1}\lambda + a_{0} + (b_{2}\lambda^2 + b_{1}\lambda + b_{0})e^{-\lambda\tau}
    \end{equation}
    where $a_{0}$, $a_{1}$, $a_{2}$, $b_{0}$, $b_{1}$, $b_{2}$ are arbitrary constants.  
    
    We call \textit{a stability switch} the change of stability of the steady state from stable to unstable or reverse. As noticed in [1], p.3928, "For nonlinear systems the occurrence of stability switch is typically associated with Hopf bifurcation, which means the appearance of periodic orbits that can be stable or unstable, depending on the type of bifurcation."
    
    In general, a characteristic function has the form $W(\lambda) = P(\lambda) + Q(\lambda)$ where $P$ is an $n$-th degree polynomial and $Q$ is a linear combination of functions $\lambda^{k} e^{-\lambda \tau_{k}}$, k=${0,1,...,n-1}$.
    
    \medskip
    \noindent
    \textbf{Lemma 1}. (application of the Mikhailov criterion). Assuming $W$ has no purely imaginary roots, the steady state of system with the characteristic equation (2) is locally stable only if 
    
    \begin{equation*}
    \Updelta_{0 \leq \omega < \infty} \arg{W(i \omega)} = \frac{\pi}{2} deg P = \frac{3 \pi}{2}
    \end{equation*}
    
    \noindent
    \textbf{Proof}. From [4] Theorem 4.5, we know that the condition sufficient and necessary for stability is that all characteristic roots have negative real parts. In other words, there should not be any roots in the right half of complex plane. 
    
    By the Argument Principle and the fact that  $W(\lambda)$ has real coefficients, the number of roots with negative real part ($N$) can found by computing $\Updelta_{C} \arg{W(i \omega)}$ on the boundary $C$ of the right half circle with center at origin and radius $R >> 0$ and multiplying the value by $-\frac{1}{2\pi}$. 
    
    By arguments similar to the proof of Rouche's Theorem and the fact that the characteristic function $W(\lambda)$ has real coefficients, the change of argument on the curvilinear part is $- \pi deg P$ and changes of two vertical half lines are equal. From this, it can be shown that the condition on roots is equivalent to the total change of argument of $W(i \omega)$ as $\omega$ increases from $0$ to $+\infty$ being $\frac{\pi}{2} deg P$. $\square$\\

    \noindent
    Stability switches can only occur if increasing delay changes the nature of roots of $W$: when $\tau=0$, $W$ has two complex conjugate roots with positive real part; the equilibrium becomes stable if 
    after the introduction of delay $W$ does not have any roots with positive real part.
    From the continuous dependence of characteristic roots on the value $\tau$, stability switches occur when the pair becomes purely imaginary. Thus a positive $\omega$ can be found such that:
    \begin{equation}
          W(i\omega) = 0 \Rightarrow P(i\omega) = -Q(i\omega) \Rightarrow |(i\omega)^3 + a_{2}(i\omega)^2 + a_{1}i\omega + a_{0}|^{2} = |b_{2} (i \omega)^2 + b_{1}i\omega + b_{0}|^{2}
    \end{equation}
    Rewriting in terms of $x = \omega^{2}$ we get a cubic polynomial $F(x)$:
    \begin{equation}
          F(x) = x^{3} +  (a_{2}^2 - 2a_{1} - b_2^2)x^{2} + (a_{1}^2 - 2a_{0}a_{2}-b_{1}^{2}-2 b_0 b_2)x + a_{0}^{2} - b_{0}^{2}
    \end{equation}
    
    \noindent
    Under some conditions which we are going to discuss below, the polynomial will have one negative root $x_{0}$ and two positive roots $x_{2}>x_{1}>0$ corresponding to $\omega_{2}>\omega_{1}$. \\

    \noindent
    \textbf{Lemma 2}.
    Let $\bar{x}$ be a positive root of $F(x)$ such that for the corresponding 
    value of $\tau$ the system gains stability. Then $F'(\bar{x})<0$, 
    the discriminant of $F(x)$
    is greater than 0 and $\lvert a_{0} \rvert > \lvert b_{0} \rvert$.\\
    
    \noindent
    \textbf{Proof}. For a stability switch to appear, we need a value of $\omega$ with $W(i\omega) =0$, hence $F(x)$ should have a corresponding positive root $\bar{x}$. As shown in [3], p.80 the sign of $F'(\bar{x})$ determines in which direction characteristic roots cross the imaginary axis in the complex plane with increasing delay. If $F'(\bar{x})<0$, the pair passes the imaginary axis from right to left (then the argument changes from $-\frac{ \pi}{2}$ to $\frac{3 \pi}{2}$), and if $F'(\bar{x})>0$, the conjugate pair moves from the left to the right part of plane (the argument decreases by $2 \pi$). Consequently, system can only gain stability if $F'(\bar{x})<0$. \\
    
    \noindent
    Since the derivative at the smallest and the greatest root of a cubic
    polynomial with leading term $x^3$ is positive, the above implies that
    $\bar{x}$ is the middle of three roots, and it is the first root to the 
    right of $0$. 
    In other words, $F(x)$ has two positive roots $x_1 = \bar{x} < x_2$, and 
    a negative root. In particular, the discriminant is positive and $F(0) > 0$, meaning that $\lvert a_{0} \rvert > \lvert b_{0} \rvert$. $\square$
    
    \begin{figure}[h!]
        \centering
        \begin{subfigure}[b]{0.18\linewidth}
            \includegraphics[width=\linewidth]{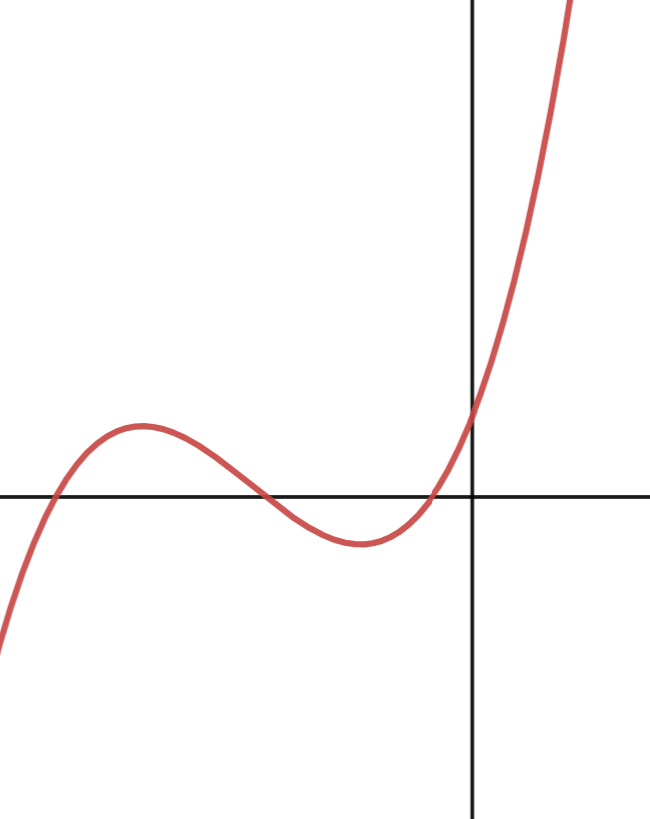}
            \caption{3 negative}
        \end{subfigure}
        \begin{subfigure}[b]{0.18\linewidth}
            \includegraphics[width=\linewidth]{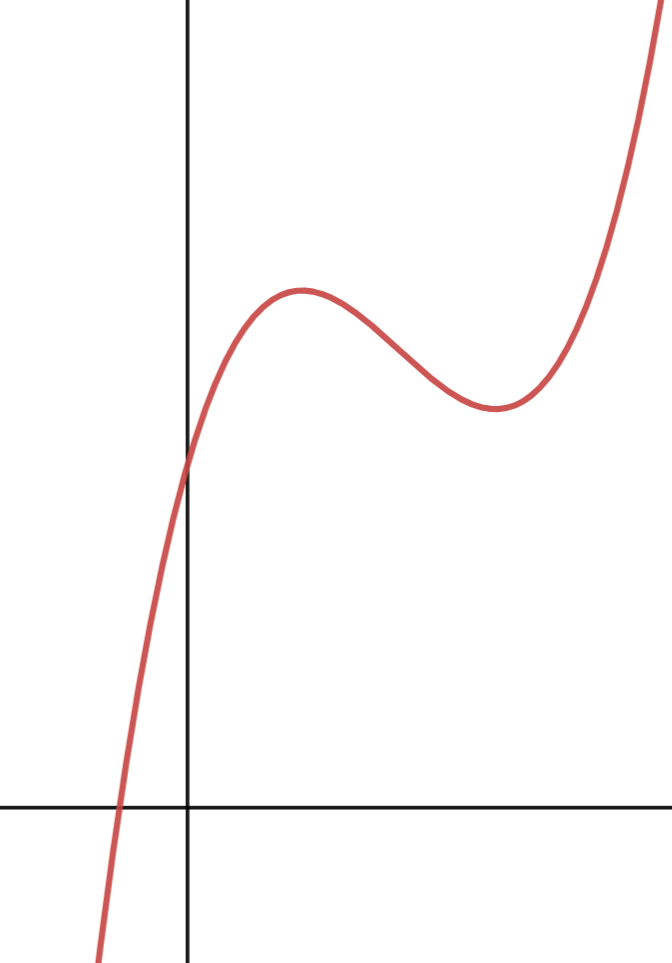}
            \caption{1 negative}
        \end{subfigure}
        \begin{subfigure}[b]{0.18\linewidth}
            \includegraphics[width=\linewidth]{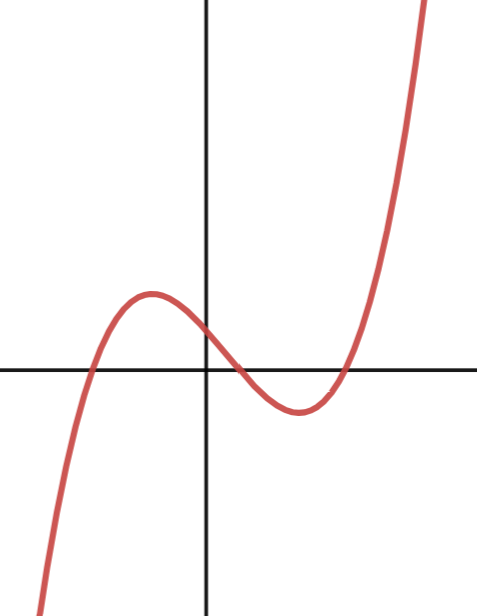}
            \caption{1 negative and \\ 2 positive}
        \end{subfigure}
        \caption{Possibilities for roots of $F(x)$.}
        \label{fig:f(x)}
    \end{figure}
    
    \medskip
    \noindent
    Lemma 3 gives simpler conditions which are  sufficient but not necessary. \\
    
    \noindent
    \textbf{Lemma 3}. If discriminant of $F(x)$
    is greater than 0 and $a_{1}^{2}-2a_{0}a_{2}-b_{1}^{2} -2 b_0 b_2 < 0$, and the Mikhailov criterion is met then stability can be gained. \\
    
    \noindent
    \textbf{Proof}. $a_{1}^{2}-2a_{2}-b_{1}^{2}-2 b_0 b_2 < 0$ means $F'(0)<0$ and positive discriminant guarantees that three roots exist. $F'(0)<0$ means there is only one negative root (this eliminates the case of three negative roots); also, as $F(x)$ has three roots, two other roots are positive, giving us two values of $\omega$ (see Fig.1c) as well as $F(0) > 0$, so $\lvert a_{0} \rvert > \lvert b_{0} \rvert$. \\
    
    \begin{figure}[h!]
        \centering
        \begin{subfigure}[b]{0.2\linewidth}
            \includegraphics[width=\linewidth, height=2cm]{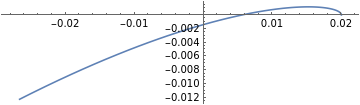}
        \end{subfigure}
        \begin{subfigure}[b]{0.1\linewidth}
            \includegraphics[width=\linewidth, height=2cm]{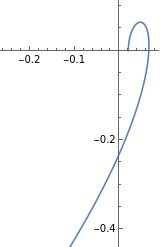}
        \end{subfigure}
        \begin{subfigure}[b]{0.18\linewidth}
            \includegraphics[width=\linewidth, height=2cm]{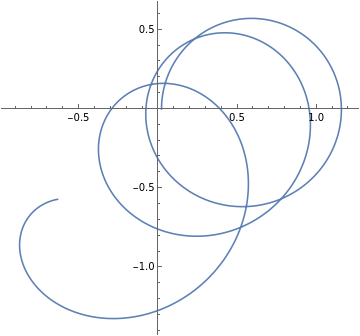}
        \end{subfigure}
        \caption{Mikhailov hodograph shows the system is unstable.}
        \label{fig:wrong mikh}
    \end{figure}
    
    \noindent
    Moreover, by adding conditions on the Mikhailov hodograph from Lemma 1 we can guarantee that stability switches will occur. As seen in Lemma 1, we want the curve traced by $W(i \omega)$ as $\omega$ increases from $0$ to $\infty$ to go around the origin counterclockwise and tend to negative part of vertical axis. By changing $\tau$ we can adjust the change of argument. At critical values of $\tau$, the hodograph should pass through the origin. For smaller values of $\tau$, the curve should cross the imaginary axis below 0. In the interval between critical values, the curve should cross the imaginary axis above 0. For bigger values of $\tau$, the curve may cross the imaginary axis below 0 or go to the right (see Fig.2). $\square$\\ 
    
    \bigskip
    \noindent
    Now we describe how to find critical values of $\tau$ if they exist. From the
    coefficients of $W(i\omega)$ we can compute the coefficients of $F(x)$ and 
    its positive roots give two values of $\omega$, hence we know the two conjugate pair of $\lambda$ that lands on the imaginary axis
    when $W(\lambda) =0$.
    By comparing real and imaginary parts $sin$ and $cos$ of $\omega \tau$ can be found.
    \begin{equation}
    \begin{aligned}
        \frac{\lambda^3 + a_{2}\lambda^2 + a_{1}\lambda + a_{0}}{-b_{0}-b_{1} \lambda-b_{2} \lambda^2} =
         e^{-i \omega \tau} =\cos{\omega \tau} - i \sin{\omega \tau} \\ 
    \end{aligned}
    \end{equation}
    
    \noindent
    Then $\omega_{j} \tau_{j} = \alpha_{j} + 2 \pi n$ with unique $\alpha_{j} \in [0; 2 \pi)$ for $n \geq 0$, $j=1,2$. So we get $\tau_{j n} = \frac{\alpha_{j}}{\omega_{j}} + \frac{2 \pi}{\omega_{j}} n$. \\
    
    \noindent
    Since $F'(x_1) < 0 < F'(x_2)$, the value $j=1$ makes characteristic roots cross the imaginary axis from right to left, and $j=2$ makes them cross in the reverse direction. Since $\omega_{1} < \omega_{2}$, the latter happens more often, so a steady state eventually becomes unstable. Consequently, the unstable steady state can initially gain stability if $\tau_{10} < \tau_{20}$ but remains unstable if $\tau_{20} < \tau_{10}$. 
    Also note that if we can find a pair $\tau_{1k}<\tau_{2k}$ with which the system gains stability, then we know that $\tau_{10} < \tau_{20}$ and initial stability interval happens on $[\tau_{10}, \tau_{20}]$.
    Therefore, existence of critical values by itself proves that the Mikhailov criterion is met. In a case when critical values cannot be easily found, we reformulate conditions from Lemma 1 and 2:
    \begin{enumerate}
        \item $F(x)$ has one negative and two positive roots; 
        \item $ \Updelta_{0 \leq \omega < \infty} \arg{W(i \omega)} = \frac{3 \pi}{2} $
    \end{enumerate}
    Splitting $W(i\omega)$ into real and imaginary parts, $W(i\omega) = W_r(\omega) + i W_i 
    (\omega)$ we introduce  assumptions that are sufficient to meet the condition 2 above:
    
    \begin{enumerate}
        \item[3.] $W_{r}(0) > 0$ and $W'_{r}(\omega) <0$ for all  $\omega > 0$;
        \item[4.] There exists such $\bar{\omega}$  that  $W_{r}(\bar{\omega}) < 0$ and $W_{i}(\omega) > 0$ for all  $\omega \in (0, \bar{\omega})$.
    \end{enumerate}

    \noindent 
    \textbf{Theorem.} Assume that for some $\tau$
    \begin{equation}
    \begin{aligned}
        0 > b_{0}> - a_{0}, \quad 0> a_{1}, \quad a_{2}> a_{1}^{2}/2 \lvert b_{0} \rvert,  
        \\
        \text{ $b_{1}$ and $b_{2}$ are sufficiently small}, 
        \\
        \text{and there exists such } \bar{\omega} \text{ that } W_{r}(\bar{\omega}) < 0
        \text{ and }W_{i}(\omega) > 0 \text{ on the interval } (0, \bar{\omega}).
    \end{aligned}
    \end{equation}
    Then stability gain occurs, after a minimal threshold $\tau_{10} \leq \sqrt{\frac{2a_{2}}{\lvert b_{0}\rvert}}$. \newline
    
    \noindent
    \textbf{Proof}. Conditions (6) imply that for $\tau=0$ the steady state is unstable: 
    the curve $W(i\omega)$ remains in Quadrants IV and III hence $\Delta\; \arg{W(i \omega)} = - \frac{\pi}{2}$ and $W(\lambda)$ has two roots in the right half plane. 
    By Lemma 3 it remains to show that the Mikhailov Criterion is met for the
    value of $\tau$ that ensures the statements about $W_i, W_r$. \\
    
    \noindent
    We are going to assume $b_{1}=0$ and $b_{2}=0$, but the proof also holds true for 
    small values (in a sense that we do not specify here), due to continuous dependence 
    of roots of the characteristic function $W$ on its parameters.
    Now we find real and imaginary parts of $W$:
    \begin{equation}
    \begin{aligned}
        W_{r}(\omega) &=  -a_{2} \omega^2 +  a_{0} + b_{0} cos(\omega \tau)\\
        W_{i}(\omega) &=  -\omega^3 +  a_{1} \omega - b_{0}  sin(\omega \tau)
    \end{aligned}
    \end{equation}
    When $\omega=0$, $W_{r}(0)=a_{0}+b_{0}>0$ and $W_{i}(0)=0$. Hence, $\arg {W(0)} = 0$.
    \begin{equation}
    \begin{aligned}
        \sin {\arg {W(\omega)}} = \frac{-\omega^3 +  a_{1} \omega - b_{0}  sin(\omega \tau)} { \sqrt{W_{r}^{2}+W_{i}^{2}}} 
        \xrightarrow{\omega \xrightarrow{} +\infty} -1 \\
        \cos {\arg {W(\omega)}}  = \frac {-a_{2} \omega^2 +  a_{0} + b_{0} cos(\omega \tau)} {\sqrt{W_{r}^{2}+W_{i}^{2}}}
        \xrightarrow{\omega \xrightarrow{} +\infty} 0
    \end{aligned}
    \end{equation}
    Therefore, $\arg {W(\omega)} \rightarrow{} \frac{3 \pi}{2}$ and the total change of argument is equal to $\frac{3 \pi}{2}+2l\pi, l \in \mathbb{Z}$. We are only left to prove that $l=0$.    Now, we show that $W_{r}$ is decreasing. 
    \begin{equation*}
        W'_{r}(\omega) =  -2 a_{2} \omega + \lvert b_{0} \rvert \tau sin(\omega \tau) \leq{} \omega (\lvert b_{0}\rvert \tau^{2}-2a_{2})<0
    \end{equation*}
    since $\sin{\omega \tau} \leq \omega \tau$ and $\tau < \sqrt{\frac{2a_{2}}{\lvert b_{0}\rvert}}$. \\
    
    \noindent
    Since we know that there is $\bar{\omega}$ such that $W_{i}(\omega)>0$ for all $\omega \in (0, \bar{\omega})$ and $W_{r}(\bar{\omega})<0$, the curve goes through Quadrants I and II and we also know that the curve tends to negative half of horizontal axis in Quadrant III; in addition to that, existence of $\bar{\omega}$ and decreasing $W_{r}$ mean that the curve doesn't intersect Quadrant IV. Thus, the change of $\arg{W(i\omega)}$ is $\frac{3 \pi}{2}$, so the steady state is stable due to the Mikhailov Criterion. The sketch of the Mikhailov hodograph is presented below. $\square$\\
    
    \begin{figure}[h!]
        \centering
        \begin{subfigure}[b]{0.45\linewidth}
            \includegraphics[width=\linewidth]{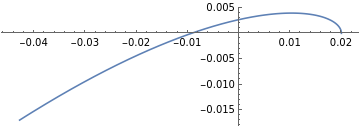}
        \end{subfigure}
        \caption{Mikhailov hodograph shows the system is stable.}
        \label{good mikh}
    \end{figure}
    
    \noindent
    We simplify the inequalities for $\bar{\omega}$ given in Theorem, approximating \textit{sin} and \textit{cos} using Taylor series. \\
    
    \noindent
    \textbf{Corollary}. Assume that
    \begin{equation*}
    \begin{aligned}
        0 > b_{0}> - a_{0}, \quad 0> a_{1}, \quad a_{2}> a_{1}^{2}/2 \lvert b_{0} \rvert
        \text{ and that }
    \end{aligned}
    \end{equation*}
    \begin{equation}
    \begin{aligned}
        g(\tau) = \frac{a_{0}+2 \lvert b_{0} \rvert}{6} \tau^{3} - \frac{\lvert a_{1} \rvert}{2} \tau^{2} - a_{2} \tau + \frac{\lvert a_{1} \rvert a_{2} + a_{0} - \lvert b_{0} \rvert}{\lvert b_{0} \rvert}
    \end{aligned}
    \end{equation}
    has at least one root in $[0, \sqrt{\frac{2 a_{2}}{\lvert b_{0} \rvert}}]$ and  $b_{1}$ and $b_{2}$ are sufficiently small.
    Then stability is gained, with the minimal thresold $\tau_{10} < \sqrt{\frac{2a_{2}}{\lvert b_{0}\rvert}}$. \newline
    
    \noindent
    \textbf{Proof}. 
    We estimate $W_{i}(\omega)$ in the following way
    \begin{equation*}
    \begin{aligned}
        W_{i}(\omega)=-\omega^3-\lvert a_{1}\rvert \omega + \lvert b_{0}\rvert  sin(\omega \tau) \geq
        -\omega^3-\lvert a_{1}\rvert \omega + \lvert b_{0}\rvert  (\omega \tau-\frac{\omega^{3} \tau^{3}}{6}),
    \end{aligned}
    \end{equation*}
    independently of $\tau$. Thus,
    \begin{equation*}
    \begin{aligned}
        W_{i}(\bar{\omega}) \geq \bar{\omega}(\lvert b_{0}\rvert \bar{\tau}-\lvert a_{1}\rvert- \bar{\omega}^2 \frac{\lvert b_{0}\rvert \bar{\tau}^{3}+6}{6}).
    \end{aligned}
    \end{equation*}
    It can be easily seen that $W_{i}(\omega)>0$ for $\omega \in (0, \bar{\omega})$ if
    \begin{equation}
    \begin{aligned}
    \bar{\omega}^{2}<\frac{6(\lvert b_{0}\rvert \bar{\tau}-\lvert a_{1}\rvert)}{\lvert b_{0}\rvert \bar{\tau}^{3}+6}
    \end{aligned}
    \end{equation}
    Next, we need $W_{r}(\bar{\omega})<0$. By approximating $cos(\bar{\omega} \bar{\tau})$ with $1 - \frac{\bar{\omega}^{2} \bar{\tau}^{2}}{2}$ we get:
    \begin{equation*}
    \begin{aligned}
    W_{r}(\bar{\omega}) =-a_{2} \bar{\omega}^2 + a_{0} - \lvert b_{0} \rvert cos(\bar{\omega} \bar{\tau}) \leq
    -a_{2} \bar{\omega}^2 + a_{0} - \lvert b_{0} \rvert + \frac{\lvert b_{0} \rvert \bar{\omega}^{2} \bar{\tau}^{2}}{2} = a_{0} - \lvert b_{0}\rvert + \omega^{2} (\frac{\lvert b_{0}\rvert \bar{\tau}^{2}}{2} - a_{2})
    \end{aligned}
    \end{equation*}
    Hence, $W_{r}(\bar{\omega})<0$ if 
    \begin{equation}
    \begin{aligned}
    \bar{\omega}^{2} > \frac{a_{0} - \lvert b_{0}\rvert}{a_{2} - \frac{\lvert b_{0}\rvert \bar{\tau}^{2}}{2}}
    \end{aligned}
    \end{equation}
    Therefore, such $\bar{\omega}$ exists if $(10) > (11)$ (notice that $(11) > 0$ so we can 'fit in' $\bar{\omega}^{2}$).
    \begin{equation*}
    \begin{aligned}
    \frac{6(\lvert b_{0}\rvert \bar{\tau}-\lvert a_{1}\rvert)}{\lvert b_{0}\rvert \bar{\tau}^{3}+6} >
    \bar{\omega}^{2} > \frac{a_{0} - \lvert b_{0}\rvert}{a_{2} - \frac{\lvert b_{0}\rvert \bar{\tau}^{2}}{2}}
    \end{aligned}
    \end{equation*}
    Then 
    \begin{equation*}
    \begin{aligned}
    g(\bar{\tau}) = \frac{a_{0}+2 \lvert b_{0} \rvert}{6} \bar{\tau}^{3} - \frac{\lvert a_{1} \rvert}{2} \bar{\tau}^{2} - a_{2} \bar{\tau} + \frac{\lvert a_{1} \rvert a_{2} + a_{0} - \lvert b_{0} \rvert}{\lvert b_{0} \rvert} < 0
    \end{aligned}
    \end{equation*}
    From $g(0)>0$, $g'(0)<0$ and positive leading coefficient we know that $g(\tau)$ has just one negative root or one negative and two positive roots. So for $g(\tau)$ to have negative values in the interval $[ 0, \sqrt{\frac{2 a_{2}}{\lvert b_{0} \rvert}}]$ one just needs to see that $g(\tau)$ has a root in  $[0, \sqrt{\frac{2 a_{2}}{\lvert b_{0} \rvert}}]$, e.g. by 
    applying Sturm's algorithm.  $\square$ \\
    
    \noindent
    \textbf{Remark.} We observe that if 
    $W_i(\omega)$ is positive on a nonempty interval $(0, \bar{\omega})$ - as it happens in the     above Theorem and Corollary - then 
    $W'_i(0) > 0$ and hence $\bar{\tau} > \frac{a_1}{b_0}$.
    
\section{Application of Theorem}
    In this section we provide an example illustrating the application of the Theorem to the system (1). For delays that are near 0, we have an unstable steady state; however, the steady state gains stability at the first threshold value of delay and loses it at the next critical value.\newline
    
    \noindent
    The characteristic function for Eqs. (1) is the following:
    \begin{equation*}
        W(\lambda) = \lambda^3 + a_{2}\lambda^2 + a_{1}\lambda + a_{0} + b_{0} e^{-\lambda\tau}
    \end{equation*}
    \begin{equation*}
    \begin{aligned}
        a_{0} = 0.16, \qquad 
        a_{1} = -0.23, \qquad 
        a_{2} = 0.97 \\
        b_{0} = -0.14 \\
        A_{1} = 1, \qquad
        A_{2} = 1, \qquad
        A_{3} = -2 \\
    \end{aligned}
    \end{equation*}
    We choose initial conditions as $x(t \leq \tau) = y(t \leq \tau) = z(t \leq \tau) = t^{2}$. 
    $F(x)$ has the following form:
    \begin{equation*}
    \begin{aligned}
        F(x) = x^{3} + 1.4009 x^2 - 0.2575x + 0.006 \\ 
    \end{aligned}
    \end{equation*}
    From roots of $F(x)$ we find values of $\omega$.
    \begin{equation*}
    \begin{aligned}
        \omega_{1} = 0.165814 \\ 
        \omega_{2} = 0.373109 \\
    \end{aligned}
    \end{equation*}
    So the critical value of $\tau$ will give $0.165814i$ and $0.373109i$ characteristic roots. 
    We find $e^{-\lambda_{1} \tau} = 0.952361-0.304972i$ and $e^{-\lambda_{2} \tau} = 0.178326-0.983971i$. By diving the value of angle in radians by $\omega_{j}$, we know critical values of $\tau$.
    \begin{equation*}
    \begin{aligned}
        \tau_{1} = 1.86902 \\
        \tau_{2} = 3.7295 \\
    \end{aligned}
    \end{equation*}

    \noindent
    Note that $1.86988 < \bar{\tau} < 3.41087$  satisfies inequalities from Corollary:
    \begin{equation*}
    \begin{aligned}
        \frac{a_{1}}{b_{0}}=1.64 < \bar{\tau}=1.87 <
        \sqrt{\frac{2a_{2}}{\lvert b_{0}\rvert}}=3.72 \\
        g(\bar{\tau}) = \frac{2431}{1400} - 0.97 \bar{\tau} - 0.115 \bar{\tau}^2 + \frac{11}{150} \bar{\tau}^3 < 0 \text{ for } \bar{\tau} \in (1.86988,3.41087)
    \end{aligned}
    \end{equation*}
    
    \noindent
    Also note that values listed above may be greater than $\tau_{j0}$; however, they guarantee existence of $\tau_{j0}$. 
    This gives the steady state $(-10.5, -50.5, 50)$. Note that in linear model appeal terms (constants) do not affect stability of the steady state but only the value of the steady state. Hence equilibrium can move from (0, 0, 0) point and 'start a love affair', but it cannot affect stability of the steady state.
    The graphs below show how emotions change with time and give phase portraits.
    When $\tau = 0$, the steady state in unstable.\\
    \includegraphics[width=5.2cm, height=3cm]{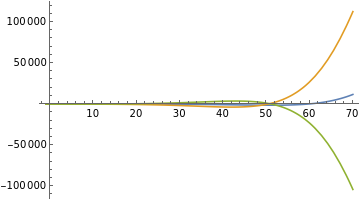}\\
    When $\tau = 1.8$, steady state is still unstable. \\
    \includegraphics[width=5.2cm, height=3cm]{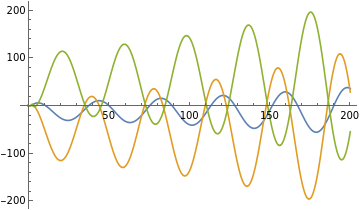}
    \includegraphics[width=5.2cm, height=3cm]{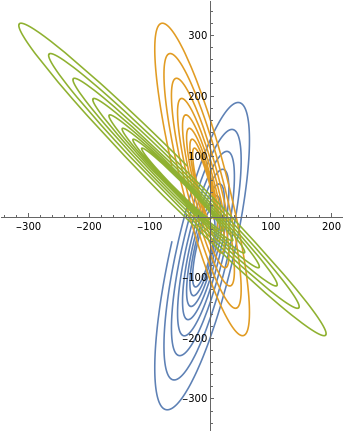} \\
    $\tau = 1.86902$ is the left boundary the system starts to gain stability. \\
    \includegraphics[width=5.2cm, height=3cm]{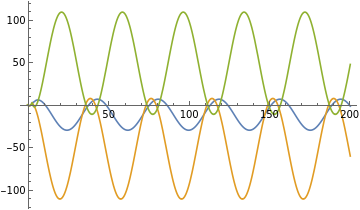}
    \includegraphics[width=5.2cm,
    height=3cm]{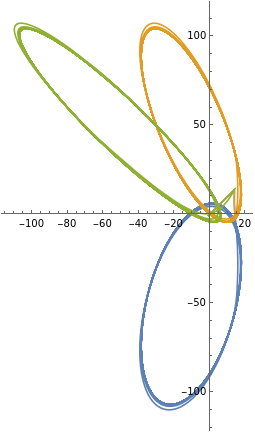} \\
    System remains stable for all $1.86902 <  \tau < 3.7295$, such as for $\tau=2.5$. \\
    \includegraphics[width=5.2cm, height=3cm]{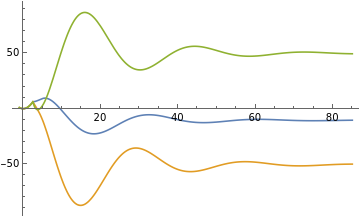}
    \includegraphics[width=5.2cm, height=3cm]{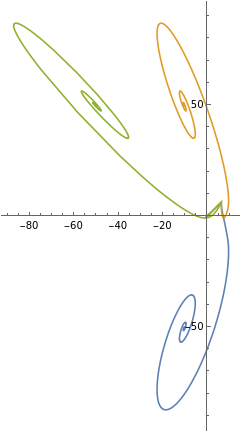} \\
    At $\tau = 3.7295$, the system loses stability. \\
    \includegraphics[width=5.2cm, height=3cm]{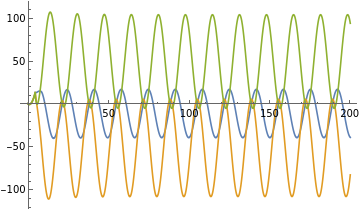}
    \includegraphics[width=5.2cm, height=3cm]{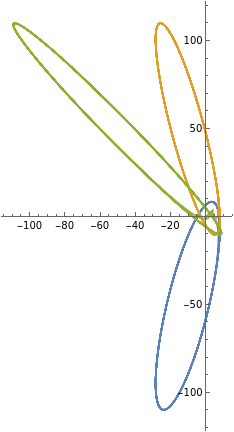}  \\
    For $\tau \geq 3.7295$, the system is unstable, such as for $\tau=4$. \\
    \includegraphics[width=5.2cm, height=3cm]{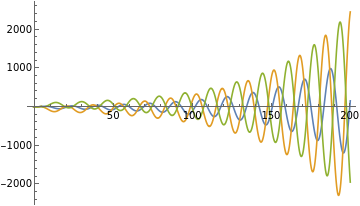}
    \includegraphics[width=5.2cm, height=3cm]{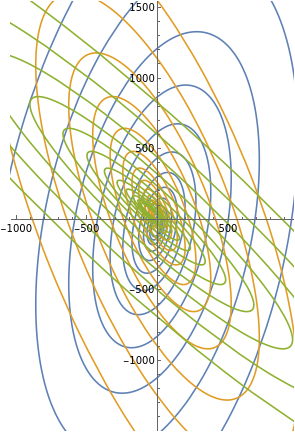}
    
\

\section{Further questions }

    In this paper we assumed $b_{2}=b_{1}=0$ and stated that the main theorem and its
    corollary hold for small enough $b_1, b_2$. It would be interesting to clarify when those values are 'small enough'. \\
    
    \noindent
    It might be interesting to look at a system with two delays, which would give a 
    more flexible model, to find how the value of one delay affects critical values of the other. Another interesting area would be behavior of the Mikhailov hodograph  for larger values of $\tau$, e.g. when it starts to spiral. 
    In this paper we wanted roots to have negative real parts; however, as can shown in Theorem 4.5 in [4], roots with zero real part do not prevent stability of the steady state as long as some
    condition is met. 
    Last but not least, systems with more agents may also be investigated: our conjecture is that four-variable case may be similar to two-variable one if agents mutually stabilize. 

\section{Appendix}
    In this section we derive the characteristic equation. We denote vector v as the following:
    \[
    v(t) = \begin{bmatrix}
        x(t)\\
        y(t) \\
        z(t)
    \end{bmatrix}
    \]
    Thus, the system in Eqs.1 can be written in terms of v.
    \[
          v'(t) = Av(t) + Bv(t-\tau) 
    \]
    \text We are looking for a special solution of the form:
    \[
    v(t) = e^{\lambda t}
    \times 
     \begin{bmatrix}
        \alpha \\
        \beta \\
        \gamma
    \end{bmatrix}
    \]
    We substitute the value of $v(t)$.
    \[
    v'(t) = \lambda e^{\lambda t} 
    \times \begin{bmatrix}
        \alpha \\
        \beta \\
        \gamma
    \end{bmatrix}
    = e^{\lambda t} A
         \times
         \begin{bmatrix}
        \alpha \\
        \beta \\
        \gamma
    \end{bmatrix}
    + e^{\lambda (t - \tau)}B
         \times
         \begin{bmatrix}
        \alpha \\
        \beta \\
        \gamma
    \end{bmatrix}
    \]
    
    \[
    (\lambda I - A - B e^{-\lambda\tau}) 
    \times
        \begin{bmatrix}
        \alpha \\
        \beta \\
        \gamma
    \end{bmatrix} = 0
    \]
    \text We assume non-trivial solution, i.e. at least one of the
    values $\alpha$, $\beta$ and $\gamma$ is not $0$.
    \\ 
    
    \noindent 
    Hence, $\det (\lambda I - A - B e^{-\lambda \tau}) = 0$.
    For the sake of simplicity, in this paper we assume that in some basis the 
    matrix $\lambda I - A - B e^{-\lambda\tau}$ has the
    form:
    \[
    M = \begin{pmatrix}
        \lambda & 0 & a_{0} + b_{0} e^{-\lambda\tau} \\
        1 & \lambda & -a_{1}-b_{1} e^{-\lambda\tau} \\
        0 & 1 & \lambda+a_{2} + b_{2} e^{-\lambda\tau}
        \end{pmatrix}
    \]
    \text From substituting $\det (M) = 0$ we get the characteristic equation of the 3-agent matrix: $ W(\lambda) = \lambda^3 + a_{2}\lambda^2 + a_{1}\lambda + a_{0} + (b_{2}\lambda^2 + b_{1}\lambda + b_{0})e^{-\lambda\tau}$. 
    In full generality the characteristic function will also have $(c_1 \lambda + c_0)e^{-2 \lambda \tau}$ and $d_0 e^{-3 \lambda \tau}$ terms.

\section{References}
    
    [1] N. Bielczyk, M. Bodnar, U. Forys, Delay can stabilize: Love affairs dynamics, Appl. Math. Comput. 219 (2012) 3923–3937.
    
    \noindent
    [2] S. Strogatz, Love affairs and differential equations, Math. Mag. 65 (1) (1988) 35.

    \noindent
    [3] K.L. Cooke, P. van den Driessche, On zeroes of some transcendental equations, Funkcj. Ekvacioj 29 (1986) 77–90.
    
    \noindent
    [4] R. Bellman and K.L. Cooke, Differential-Difference equations, R-374-PR (1963)
    
    \noindent
    [5] T. Erneux, Applied Delay Differential Equations (2009)
   
\end{document}